\documentclass[11pt]{amsart}
\usepackage{geometry,verbatim,xcolor}                
\usepackage{graphicx}
\usepackage{amssymb}
\usepackage{epstopdf}
\usepackage{lscape}
\newtheorem{example}{Example}

\def\R{\Bbb R}

\def\spp{\vspace{5pt}\noindent}
\newtheorem{lemma}{Lemma}

\newtheorem{theorem}{Theorem}
\newtheorem{corollary}{Corollary}

\title{Second Order Spiral Splines}
\author{Lyle Noakes}
\date{January 10, 2019}                                           
\address{Lyle.Noakes@uwa.edu.au (Department of Mathematics \& Statistics, The University of Western Australia, 35 Stirling Highway, Crawley, WA 6009, AUSTRALIA)}

\begin{document}

\begin{abstract}
Second order spiral splines are $C^2$ unit-speed planar curves that can be used to interpolate a list $Y$ of $n+1$ points in $\R ^2$ at times specified in some list $T$, where $n\geq 2$. Asymptotic methods are used to develop a fast algorithm, based on a pair of tridiagonal linear systems and standard software. The algorithm constructs a second order spiral spline interpolant for data that is convex and sufficiently finely sampled. \end{abstract}

\maketitle


\section{Introduction}\label{intro}
Given a finite sequence $T$ of real numbers $T_0<T_1<\ldots <T_n$, a pair $V$ of vectors $V_0,V_n\in \R ^m$ and a finite sequence $Y$ of points 
$Y_0,Y_1,\ldots ,Y_n\in \R ^m$ where $m\geq 1$, an {\em interpolant} of $(Y,V)$ at $T$ is defined to be a $C^2$ curve $y:[T_0,T_n]\rightarrow \R ^m$ satisfying $y(T_j)=Y_j$ for $0\leq j\leq n$, as well as the auxiliary conditions  
$y^{(1)}(T_0)=V_0$ and $y^{(1)}(T_n)=V_n$ on derivatives $y^{(1)}$ of $y$ at $T_0$ and $T_n$. Then $y$ is said to  {\em interpolate} $(Y,V)$ at $T$. 

\spp
A standard interpolant is the minimiser $y$ of   
$$J(y)~:=~\int _{T_0}^{T_n}\Vert y^{(2)}(t)\Vert ^2~dt$$
where $\Vert ~\Vert$ denotes the Euclidean norm, and $y^{(2)}$ is the second derivative. There is precisely one minimiser, namely the {\em cubic polynomial spline} which is the $C^2$ piecewise-cubic polynomial interpolant with knots at $T$.  
The cubic polynomial spline is easily calculated from $T$ and $(Y,V)$, by solving a tridiagonal system of linear equations \cite{deBoor}. It is much more difficult, both from a theoretical and practical point of view,  to find $C^2$ interpolants that are  {\em unit-speed}, namely $\Vert y^{(1)}(t)\Vert =1$ for all $t\in [T_0,T_n]$. 

\spp
Unit-speed curves are significant for
highway and railway design \cite{meek0}, \cite{meek}, motion planning
for robots \cite{Kelly}, and path planning for unmanned aerial vehicles  and military aircraft \cite{dai}, \cite{looker}.  Whereas cubic polynomial splines always exist, an evident necessary condition for existence of a $C^2$ unit-speed interpolant is 
\begin{equation}\label{nc}\Vert Y_j-Y_{j-1}\Vert \leq L_j:=T_j-T_{j-1}\quad \hbox{for ~}1\leq j\leq n. \end{equation}

\spp
By an {\em elastic spline} we mean\footnote{The term {\em elastic spline} is sometimes used differently, to mean an interpolant where $T$ is not specified in advance \cite{birkhoff}.} a critical point of the restriction of $J$ to the space of unit-speed interpolants. An elastic spline is known to be precisely a $C^2$ interpolant that is a track-sum $y$ of elastica (elastic curves) $y_j:[T_{j-1},T_j]\rightarrow \R ^m$. An {\em elastica}  is defined as\footnote{The length of $y_j$ is then $T_j-T_{j-1}$. Sometimes the term elastica has a different meaning, where the length of $y_j$ is not prescribed in advance \cite{birkhoff}. It seems that the plural of elastica is elastica.} a critical point of $J$ restricted to unit-speed curves $y_j$ 
where $y_j(T_{j-1}),~y_j^{(1)}(T_{j-1}),~y_j(T_{j}),~y_j^{(1)}(T_{j})$ are all prescribed in advance. 
Elastic splines are extensively studied, with well-developed algorithms for interpolation, as discussed in \S\ref{backsec}.  These algorithms for elastic splines require solutions of systems of nonlinear equations for variable parameters, and consequently very significant effort, compared with the sparse linear systems needed for interpolation by (non-unit-speed) polynomial splines.   

\spp
Besides elastica there are simpler classes of unit-speed track-summands called {\em polynomial spirals} \cite{Kelly}.  The well-studied class of {\em clothoidal splines},  whose curvatures are $C^0$ and piecewise-affine, is insufficiently rich for interpolation when $T$ is given in advance.  In the present paper we suppose that the inequalities (\ref{nc}) hold strictly, and consider the larger class of {\em second order spiral splines}, namely interpolants that are $C^2$ track-sums $y_\theta$ of {\em second order generalised Cornu spirals}, whose curvatures are $C^0$ and {\em piecewise-quadratic} in the parameter $t$. 

\spp
In the present paper the data $(Y,V)$ is generated by unknown strictly convex\footnote{The nonconvex case is more complicated, with potentially $2^n$ interpolants.} curves $x$ with suitably bounded derivatives, and sufficiently fine sampling specified by $T$. The setting is described in detail in \S \ref{admiss}.  
In \S \ref{newsec}, \S \ref{appsec} a pair of tridiagonal linear systems defines a $C^1$ cubic polynomial spline $\hat \theta :[t_0,T_n]\rightarrow \R $.  Our construction echoes the well-known algorithm for polynomial cubic splines, where there is a single tridiagonal system \cite{deBoor}. 
As with the classical algorithm, our linear systems are well-conditioned and very quick to solve. 
Our construction is more than twice as complicated as the classical algorithm, and requires more care to implement, but the description in \S \ref{newsec}, \S \ref{appsec} is straightforward. 

\spp
Unfortunately the same cannot be said for the proof of  Theorem \ref{thm0}, which is our main result. 
The proof, spanning \S \ref{asympsec}, \S \ref{newnewsec}, \S \ref{newnewsec1}, \S \ref{completepfsec} uses a number of difficult asymptotic expansions, followed by a sequence of careful manipulations and many equations. 
Theorem \ref{thm0} says that, relative to a parameter $\epsilon >0$ measuring distances between successive data, we have $\hat \theta =\theta +O(\epsilon ^4)$ where $y_\theta$ is the desired interpolant. This immediately gives rise to some {\em approximate interpolants} by spiral splines, as illustrated in Examples \ref{ex1}, \ref{ex2}, \ref{ex3}. More significantly, as demonstrated later in \S \ref{gapsec}, availability of a suitable initial guess is absolutely key to the nonlinear problem of determining $\theta$ and thereby $y_\theta$.  The principle is illustrated in  Examples \ref{ex4}, \ref{ex5}, \ref{ex6}. 
\section{More Background on Unit-Speed Interpolants}\label{backsec}
A natural class of unit-speed interpolants is that of the unit-speed {\em reparameterised polynomial splines}, studied in \cite{WKA}, \cite{eberly}, \cite{peterson} where there are algorithms for efficient computation. 
Whereas a polynomial spline can be found to interpolate $(Y,V)$ at $T$, the unit-speed reparameterisation interpolates $(Y,V)$ at a different set of parameter values. This does not work well when $T$ is given in advance, as in our situation. 
Similar difficulties arise with {\em clothoidal splines}\footnote{A clothoidal spline is a unit-speed $C^2$ planar curve $t\mapsto y(t)$ whose curvature is $C^0$ and piecewise-affine in $t$.} \cite{stoer}, \cite{bertails}, \cite{baran}, \cite{meek0}, \cite{coope} and their $C^1$ generalisations \cite{looker}. The algorithms in these papers are of significant interest, but do not apply when $T$ is already prescribed. 

\spp
Interpolation by unit-speed curves can also be performed using {\em elastic splines} \cite{leefor}, \cite{edwards}, namely\footnote{Elastic splines and elastica mean different things depending on the context \cite{linners}. In the present paper, elastic splines are {\em pinned} and {\em clamped}, and elastica are {\em fixed-length}. In \cite{leefor}, $T$ is not given in advance and the elastic splines are {\em sliding}. } critical points of the restriction of the functional $J$ to the space of {\em unit-speed} $C^2$ interpolants. An elastic spline turns out to be a $C^2$ track-sum $y$ of {\em elastica} $y_j:[T_{j-1},T_j]\rightarrow \R ^m$, satisfying the auxiliary end conditions. Elastica (elastic curves) have a long and interesting history \cite{levien}, \cite{birkhoff}, \cite{dan}, and their study reduces\footnote{On the other hand, the study of elastic splines does not reduce in this way.} to the case where $m=3$. The elastica are completely known in terms of elliptic functions \cite{singer}, with simplifications for $m=2$, which is the case of interest for the present paper. 
Although elastic splines (sometimes called {\em nonlinear splines} or {\em true splines}) are highly regarded\footnote{The objection raised at the end of \cite{leefor} does not really apply when $T$ is given.} as interpolants, they are less widely used than cubic polynomial splines. This is because the interpolation conditions for elastic splines require the solution of a system of nonlinear equations that is even more complicated\footnote{As noted on p.184 of \cite{edwards}, clothoidal splines are sometimes used to construct initial guesses for the computation of elastic splines.} than for clothoidal splines. However, the relative sophistication of elastica 
compared with cubic polynomials is not the main difficulty in computing elastic splines. As with most nonlinear problems, indeed also with our present task of finding spiral spline interpolants, the main difficulty is the construction of a suitable initial guess. 

\spp
The algorithms in \cite{leefor}, \cite{edwards} and \cite{jerome} are for {\em sliding elastic splines}, where $T$ is not given in advance.  Another condition that is sometimes imposed is that the total time $T_n-T_0$ should be fixed or, equivalently, the length of the spline is prescribed in advance. In the present paper the entire sequence $T$ is assumed to be given in advance, by analogy with standard interpolation by cubic polynomial splines. 

\spp
What is different from cubic polynomial splines is that our interpolants are required to be unit-speed. Because of this additional requirement, once $T_0$ is given, fixing the rest of the sequence $T$ is equivalent to prescribing the lengths $L_j$ of the interpolant between all consecutive data points. In these circumstances, the main contribution to interpolation by elastic splines is Theorem 4.1 of \cite{linners}, where a nonlinear system of equations is given for determining elastic splines in great generality.  

\spp
As usual with nonlinear problems, an initial guess is required in order to get started. Unfortunately \cite{linners} says very little about how this might be constructed. The present paper contributes something in this respect, although our focus is primarily on spiral splines, whereas \cite{linners} is concerned with elastic splines.  Our guess hinges on the construction of $\hat \theta$ in \S \ref{newsec}, \S \ref{appsec}. In \S \ref{gapsec} standard software is used to improve $\hat \theta$ to $\theta$. 

\spp 
We also mention the {\em discrete elastic splines} of \cite{bruckstein}, where a discrete analogue of $J$ is optimised with respect to a variable finite sequence of points approximating $y$. In effect optimisation  of $J$ with respect to $y$ is replaced by a large finite-dimensional optimisation, whose outcome depends critically on an unspecified initial guess in a high dimensional space.

\spp
In summary, reparameterised polynomial splines and clothoidal splines are unsuitable for interpolating $(Y,V)$ by $T$. Elastic splines come with non-negligible computational issues, especially the generation of suitable initial guesses.

\section{Convex Generators and Admissible Samplings}\label{admiss}
Let $n\geq 2$ be given, together with $T_0<T_n$, $0<C<B_2$, and $B_i$ where $3\leq i\leq 5$. A unit speed $C^\infty$ curve $x:[T_0,T_n]\rightarrow \R ^2$ is called a {\em convex generator} when $C<\Vert x^{(2)}\Vert _{\infty}\leq B_2$ and  
$\Vert x^{(i)}\Vert _\infty \leq B_i$ for $i=3,4,5$. For brevity we write $X$ for the space $X_{T_0,T_n,C,B_2,B_3,B_4,B_5}$ of all convex generators. 

\spp 
Let $0<A_m<A_M$ be given, together with $\epsilon >0$ chosen small in comparison with the constants $C,B_3,B_4,B_5$ that define $X$. We define  
$\mathcal{T}_\epsilon $ to be the set of 
finite sequences $T$ of reals $T_0<T_1<T_2<\ldots <T_n\leq b$ such that 
$A_m\epsilon <L_j<A_M\epsilon$ for all $1\leq j\leq n$. 

\spp
A pair  $(Y,V)$ is called {\em admissible}\footnote{The convexity assumption built into admissibility is the reason that we are able to find a unique estimate $\hat \theta$ (otherwise there would be many).} with respect to $T\in \mathcal{T}_\epsilon$ when 
$Y=(x(T_0),x(T_1),\ldots ,x(T_n))$ and $V=(x^{(1)}(T_0),x^{(1})(T_n)$ for some convex generator $x\in X$.   
Then, for $1\leq j\leq n$, set  $q_j:=(Y_j-Y_{j-1})/L_j$, $r_j:=\Vert q_j\Vert$, and 
$$k_j~:=~\frac{\sqrt{12(1-r_j^2)}}{L_jr_j}.$$
\begin{lemma}\label{bdsklem}  For small $\epsilon >0$, $r_j$ is bounded away from $0$, and $k_j\geq C+O(\epsilon )$ for all $1\leq j\leq n$. 
\end{lemma}

\spp
{\bf Proof:} We have  
~$\displaystyle{x(T_j)-x(T_{j-1})=x^{(1)}(T_{j-1})L_j+x^{(2)}(T_{j-1})\frac{L_j^2}{2}+x^{(3)}(T_{j-1})\frac{L_j^3}{6}+O(\epsilon ^4)}$ ~  where $x\in X$. In particular $L_jr_j=1+O(\epsilon )$, and $r_j$ is bounded away from $0$ when $\epsilon >0$ is small. 

\spp
Because $\Vert x^{(1)}\Vert =1$, we also have $\langle x^{(1)},x^{(2)}\rangle =0$ and $\langle x^{(1)},x^{(3)}\rangle =-\Vert x^{(2)}\Vert ^2$, where $\langle ~,~\rangle$ is the Euclidean inner product.  
Therefore, and by convexity,   
$$r_j^2L_j^2=L_j^2-\Vert x^{(2)}(T_{j-1})\Vert ^2\frac{L_j^4}{12}+O(\epsilon ^5)~\Longrightarrow ~1-r_j^2=\Vert x^{(2)}(T_{j-1})\Vert ^2\frac{L_j^2}{12}+O(\epsilon ^3)\geq \frac{C^2L_j^2}{12}+O(\epsilon ^3).$$
So for small $\epsilon $, $0<r_j<1$ and $k_j^2>12(1-r_j^2)/L_j^2\geq C^2+O(\epsilon)$.  
\qed  

\spp
From admissibility\footnote{ Looking forward, admissibility is needed to guarantee the conclusions of Theorem \ref{thm0}, but  
may be troublesome to check in applications. In such   
cases (including the examples of the present paper) Theorem \ref{thm0} may be taken to assert that when its conclusions are false, the sampling of $x$ is too irregular or too sparse. The cure, which seems to  be rarely needed, is to sample more regularly and more often from $x$.}  the $Y_j$ are separated by no more than $A_M\epsilon$, and lie in the closed disc of radius $nA_M\epsilon$ centred on $Y_0$.

\spp
Assuming $\epsilon $ is sufficiently small for the conclusion of Lemma \ref{bdsklem} to hold,  write $q_j=r_j(\cos \omega _j,\sin \omega _j)$. We also require $\omega _1\in (-\pi ,\pi )$, and $\omega _j=\omega _{j-1}+O(\epsilon)$ for $2\leq j\leq n$, so that the $\omega _j$ are uniquely determined. By admissibility, 
$V_0=(\cos \nu _0,\sin \nu _0)$ and $V_n=(\cos \nu _n,\sin \nu _n)$
where $\nu _0=\omega _1+O(\epsilon)$ and $\nu _n=\omega _n+O(\epsilon )$.

\spp
Let $\theta :[T_0,T_n]\rightarrow \R$ be a $C^1$ cubic spline with knots $T_j$ such that $\theta (T_0)=\nu _0$ and $\theta (T_n)=\nu _n$. For $1\leq j\leq n$ and all $s\in [0,L_j]$ write  $\theta _j(s):=\theta (T_{j-1}+s)=a_j+b_js+c_js^2+d_js^3$. 
Besides the two auxiliary end conditions, 
\begin{equation}\label{aux}a_1=\nu _0,~\quad \quad ~a_n+b_nL_n+c_nL_n^2+d_nL_n^3=\nu _n,\end{equation}
there are $2n-2$ conditions for $\theta$ to be $C^1$, namely 
\begin{eqnarray}
\label{eqct}a_{j+1}&=&a_j+b_jL_j+c_jL_j^2+d_jL_j^3,\\
\label{eqc1}b_{j+1}&=&b_j+2c_jL_j+3d_jL_j^2,
\end{eqnarray} 
where $1\leq j\leq n-1$. So we have $2n$ affine equality constraints on $4n$ coefficients $a_j,b_j,c_j,d_j$.   

\spp
The $C^1$ spline $\theta $ determines a {\em second order spiral spline} $y_\theta :[T_0,T_n]\rightarrow \R ^2$ given by  
\begin{equation}\label{ydf}y_\theta (t)~:=~Y_0+\int _{T_0}^t(\cos \theta (s),\sin \theta (s))~ds.\end{equation}
We see that $y_\theta$ is unit-speed with $y^{(1)}(T_0)=V_0$ and $ y^{(1)}(T_n)=V_n$. We are interested in finding $\theta$ such that  $y_\theta$ is an interpolant of $(Y,V)$ at $T$, namely  
 for $1\leq j\leq n$, 
\begin{equation}\label{placeeq}
\int _0^{L_j}(\cos \theta _j(s),\sin \theta _j(s))~ds~=~Y_j-Y_{j-1}.
\end{equation}

\spp
Condition (\ref{placeeq}) amounts to  $2n$ non-affine conditions on the $4n$ coefficients, making $4n$ conditions in total. 
\section{Tridiagonal Systems}\label{newsec}
Given $(Y,V)$ admissible with respect to $T\in {\mathcal T}_\epsilon $ where $\epsilon >0$ is small, form the $n\times n$ tridiagonal matrix 
$$
{\bf \hat T}:=\left[ \begin{array}{ccccccccc}
3L_1&-L_1&0&0&0&0&0&\ldots &0\\
-L_1&3(L_1+L_2)&-L_2&0&0&0&0&\ldots &0\\
0&-L_2&3(L_2+L_3)&-L_3&0&0&0&\ldots &0\\
0&0&-L_3&3(L_3+L_4)&-L_4&0&0&\ldots &0\\
\vdots&\vdots&\vdots&\vdots&\vdots&\vdots&\vdots&\vdots\\
\vdots&\vdots&\vdots&\vdots&\vdots&\vdots&\vdots&\vdots\\
0&0&0&0&0&0&-L_{n-2}&3(L_{n-2}+L_{n-1})&-L_{n-1}\\
0&0&0&0&0&0&0&-3L_{n-1}&9L_{n-1}+8L_n
\end{array}\right] $$
and, if $n\geq 3$, the tridiagonal matrix 
$${\bf \tilde T}:=\left[ \begin{array}{ccccccccc}
2L_1&L_1&0&0&0&0&0&\ldots &0\\
L_1&2(L_1+L_2)&L_2&0&0&0&0&\ldots &0\\
0&L_2&2(L_2+L_3)&L_3&0&0&0&\ldots &0\\
0&0&L_3&2(L_3+L_4)&L_4&0&0&\ldots &0\\
\vdots&\vdots&\vdots&\vdots&\vdots&\vdots&\vdots&\vdots\\
\vdots&\vdots&\vdots&\vdots&\vdots&\vdots&\vdots&\vdots\\
0&0&0&0&0&0&L_{n-2}&2(L_{n-2}+L_{n-1})&L_{n-1}\\
0&0&0&0&0&0&0&L_{n-1}&2L_{n-1}+3L_n/2
\end{array}\right] .$$
Solve
\begin{equation}\label{Ttileq} {\bf \tilde T}~\tilde b ~=~\tilde {\bf R}~:=~6\left[ \begin{array}{c}\omega _1-\nu _0\\
\omega _2-\omega _1\\
\omega _3-\omega _2\\
\omega _4-\omega _3\\
\vdots\\
\vdots\\
\omega _{n-1}-\omega _{n-2}\\
(3\omega _n-\nu _n)/2-\omega _{n-1}
\end{array}\right] \end{equation}
for an $n$-dimensional column vector $\tilde b$. 

\spp
Taking $\sigma$ to be the sign of the determinant of the $2\times 2$ matrix  $[ V_0~\vdots ~Y_1-Y_0]$, 
set    
\begin{eqnarray}
\label{rhojdf} \rho _j&:=&\sigma k_j\sqrt{1-\frac{k_j^2L_j^2}{20}-\frac{(\tilde b_{j+1}-\tilde b_j)^2}{60k_j^2}}\quad \hbox{for }1\leq j\leq n-1,\quad \hbox{and}\\
\label{rhondf} \rho _n&:=&\sigma k_n\sqrt{1-\frac{k_n^2L_n^2}{20}-\frac{9}{60k_n^2}\left( \frac{\omega _n-\nu _n}{L_n}+\frac{\tilde b_n}{2} \right) ^2}.
\end{eqnarray}
Then solve 
\begin{equation}\label{Thateq} {\bf \hat T}~\hat b~=~\left[ \begin{array}{c}24(\omega _1-\nu _0)-10L_1\rho _1\\
24(\omega _2-\omega _1)-10(L_1\rho _1+L_2\rho _2)\\
24(\omega _3-\omega _2)-10(L_2\rho _2+L_3\rho _3)\\
24(\omega _4-\omega _3)-10(L_3\rho _3+L_4\rho _4)\\
\vdots\\
\vdots\\
\quad \quad 24(\omega _{n-1}-\omega _{n-2})-10(L_{n-2}\rho _{n-2}+L_{n-1}\rho _{n-1})\\
\quad 24(2\omega _n+\nu _n -3\omega _{n-1})-10(3L_{n-1}\rho _{n-1}+4L_n\rho _n)
\end{array}\right] \end{equation}
for an $n$-dimensional column vector $\hat b$. 
Set
\begin{eqnarray}
\label{bcjeq}\hat c_j&:=&\frac{-7\hat b_j-3\hat b_{j+1}}{4L_j}~+\frac{5\rho _j}{2L_j}\quad \hbox{for }1\leq j\leq n-1,\\
\label{bdjeq}\hat d_j&:=&
\frac{5(\hat b_j+\hat b_{j+1})}{6L_j^2}\quad -\frac{5\rho _j}{3L_j^2}\quad \hbox{for }1\leq j\leq n-1,\\
\label{bcneq}\hat c_n&:=& \frac{6(\omega _n-\nu _n)}{L_n^2}-\frac{2\hat b_n}{L_n}+\frac{5\rho _n}{L_n},\\
\label{bdneq}\hat d_n&:=&-\frac{20(\omega _n-\nu _n)}{3L_n^3}+\frac{10\hat b_n}{9L_n^2}-\frac{40\rho _n}{9L_n^2},  
\end{eqnarray}
and finally $\hat a_1:=\nu _0$ with 
\begin{equation}\label{bajeq} \hat a_{j+1}~:=~\hat a_j+\hat b_jL_j+\hat c_jL_j^2+\hat d_jL_j^3\quad \hbox{for }1\leq j\leq n-1.\end{equation}

\spp
Then, for  all $1\leq j\leq n$, define $\hat \theta _{j}:[0,L_j]\rightarrow \R$ by  
~$\hat \theta _{j}(s):=\hat a_j+\hat b_js+\hat c_js^2+\hat d_js^3$. A $C^1$ cubic polynomial spline $\hat \theta  :[T_0,T_n]\rightarrow \R$ is then given by $\hat \theta  (T_n):=\hat \theta _{n}(L_n)$, and  
$\hat \theta  (t)~:=~\hat \theta _{j}(t-T_{j-1})$ 
for $t\in [T_{j-1},T_j)$. 
\begin{theorem}\label{thm0} $\hat \theta (T_0)=\nu _0,\hat \theta (T_n)=\nu _n$ and, for small $\epsilon >0$, 
$\Vert \theta -\hat \theta \Vert _\infty =O(\epsilon ^4)$. 
\end{theorem}

\spp
\begin{corollary}\label{cor0} $y_{\hat \theta}$ satisfies the auxiliary conditions (\ref{aux}), and $\Vert y_\theta -y_{\hat \theta }\Vert _\infty =O(\epsilon ^5)$. \qed
\end{corollary}

\spp
To prove Theorem \ref{thm0} it suffices to show that, for  all $1\leq j\leq n$, 
\begin{eqnarray}
\label{ajhat}a_j&=&\hat a_j+O(\epsilon ^4)\\
\label{bjhat}b_j&=&\hat b_j+O(\epsilon ^3)\\
\label{cjhat}c_j&=&\hat c_j+O(\epsilon ^2)\\
\label{djhat}d_j&=&\hat d_j+O(\epsilon ).
\end{eqnarray}
We now prove some asymptotic lemmas, for use in \S \ref{newnewsec}, \S \ref{newnewsec1} where (\ref{bjhat}) is proved in two steps.  In \S \ref{completepfsec}  it turns out that (\ref{bjhat}) is the key ingredient for an otherwise easy proof of  (\ref{ajhat}), (\ref{cjhat}), (\ref{djhat}), and this completes the proof of Theorem \ref{thm0}.

\section{Approximation and Interpolation}\label{appsec}
By Corollary \ref{cor0}, the $C^2$ second order spiral spline $y_{\hat \theta}$ satisfies $y_{\hat \theta}(T_0)=Y_0$, 
$y_{\hat \theta}^{(1)}(T_0)=V_0$, $y_{\hat \theta}^{(1)}(T_n)=V_n$ and $y_{\hat \theta})(T_j)=Y_j+O(\epsilon ^5)$ for $1\leq j\leq n$. Usually $y_{\hat \theta}(T_j)\not= Y_j$ for $1\leq j\leq n$. So $y_{\hat \theta}$ is only {\em approximately} an interpolant of $(Y,V)$. 

\spp
Because of the way $y_\theta$ is constructed from $\theta$, interpolation errors at $T_j$ tend to accumulate as $j$ increases. This can be avoided by using $\theta$ differently to define  $\bar y_{\theta}:[T_0,T_j]\rightarrow \R ^2$ by $\bar y_{\theta}(T_0)=Y_0$ and for $1\leq j\leq n$,  
$$\bar y_{\theta}(t)~:=~Y_{j-1}+\int _{T_{j-1}}^t(\cos \theta (s),\sin \theta (s))~ds\quad \hbox{for }T_{j-1}<t<T_j.$$

\spp
By construction $\bar y_{\theta}(T_0)=Y_0$, and $\bar y_\theta (T_j^-)=Y_j$ for all $j$, regardless of $\theta$. However unless $y_\theta$ interpolates $(Y,V)$ exactly at $T$ the left-interpolant $\bar y_\theta$ is not 
even a second order spiral spline, because $\bar y_\theta$ is only left-continuous but not continuous at $T_j$ for $1\leq j\leq n-1$.

\spp
By Theorem \ref{thm0},  $\bar y_{\hat \theta}(T_j)=Y_j+O(\epsilon ^5)$ for $1\leq j\leq n$, and 
$\bar y_{\hat \theta}^{(1)}(T_0)=V_0$, $\bar y_{\hat \theta}^{(1)}(T_n)=V_n$. So the left-interpolant $\bar y_{\hat \theta}$ 
{\em approximately} interpolates $(Y,V)$ at $T$. 
Because the errors in $\bar y_{\hat \theta}(T_j)\approx Y_j$ do not accumulate as $j$ increases,  $\bar y_{\hat \theta}$ is better than $y_{\hat \theta}$ for diagnostics.
\begin{example}\label{ex1} Set $n=10$ with $T=(0, 0.575, 0.92, 1.265, 1.38, 1.61, 1.794, 1.955, 2.07, 2.185, 2.3)$, and let the convex generator be the second order spiral $x:[0,2.185]\rightarrow \R ^2$ given by $x(0)=(0,0)$ and 
$x^{(1)}(t)=(\cos \psi (t),\sin \psi (t))$ with $\psi (t)=t+t^2+t^3$. This generates data\footnote{We have chosen $T$ in order to sample from $x$ somewhat more frequently where curvature is large. Finer sampling always helps. } $(Y,V)$ where 
$$Y~=~((0, 0), (0.494229, 0.234572), (0.436814, 0.544934), 
(0.196378, 0.419176), (0.26536, 0.333427), (0.387233, 0.460597), 
(0.245648, 0.44385),$$
$$ (0.348732, 0.376214), (0.340667, 0.473069), 
(0.262317, 0.421501), (0.348435, 0.395382))$$ 
and $V=((1,0),),(0.615762, 0.787932))$. 
%
Computation of the $4n$ coefficients for $\hat \theta $ took 0.000918 seconds in Mathematica on a 2015 a 2.2GHZ MacBook Air with 8 GB RAM. The plot of $y_{\hat \theta}$ is shown in Figure \ref{fig2}, together with the data.
\begin{figure}[htbp] 
   \centering
   \includegraphics[width=7in]{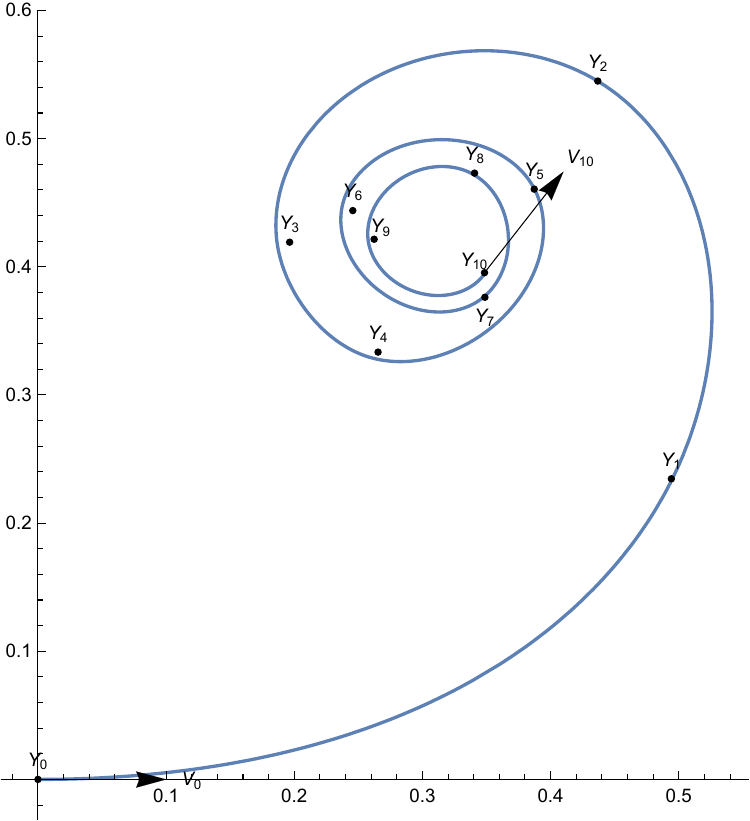} 
   \caption{$y_{\hat \theta}$ for Example \ref{ex1}}
   \label{fig2}
\end{figure} 
Although $\psi$ is just a cubic polynomial, the $C^1$ cubic spline $\hat \theta$ is not exactly $\psi$, as can be seen from the failures evident in Figure \ref{fig2} of $y_{\hat \theta}$ to interpolate at $T_3,T_4,T_6,T_9$. 
However $y_{\hat \theta}$ is not a bad initial guess for an interpolant. Finer sampling would improve it. Coarser sampling would make it worse. 
$\qed$
\end{example}  

\spp
In Example \ref{ex1} the convex generator $x$ is actually a second order spiral. Similar results are obtained when the data is generated by a more complicated convex curve. 
\begin{example}\label{ex2} Take $T=(0, 0.6, 0.96, 1.32, 1.44, 1.68, 1.872, 2.04, 2.16, 2.28, 2.4)$ and 
let $x:[0,2.4]\rightarrow \R ^2$ be given by $x(0)=(0,0)$ and $x^{(1)}(t)=e^t+t (1 + \sin (5t)/2)$. Using $x$ to generate data $(Y,V)$, Mathematica takes 0.000998 seconds to compute $\hat\theta$.  The plot of $y_{\hat \theta}$ is shown in Figure \ref{fig3}, together with the data.
%
%
\begin{figure}[htbp] 
   \centering
   \includegraphics[width=7in]{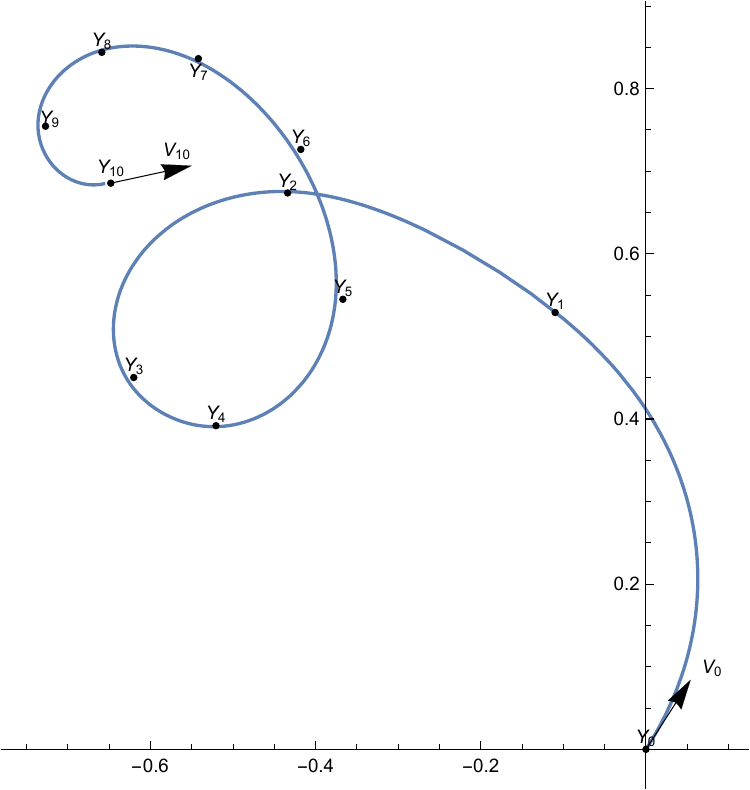} 
   \caption{$y_{\hat \theta}$ for Example \ref{ex2}}
   \label{fig3}
\end{figure}
In Figure \ref{fig3}, $y_{\hat \theta}$ almost interpolates the data, except near $T_3$, $T_5$, $T_6$, $T_9$. \qed
\end{example}
\begin{example}\label{ex3} Whereas in Example \ref{ex2} $y_{\hat \theta}$ is nearly acceptable, the quality of the approximate interpolant falls away very quickly if we use define $(Y,V)$ in the same way but over a larger interval $[T_0,T_n]$. Here we take $T=(0, 0.675, 1.08, 1.485, 1.62, 1.89, 2.106, 2.295, 2.43, 2.565, 2.7)$, and $x$ is given by the same formula over the larger interval $[0,2.7] $. Computation of $\hat \theta$ took $0.000986$ seconds, and the plot of $y_{\hat \theta}$ is shown in Figure \ref{fig4}, together with the data.
%
%
\begin{figure}[htbp] 
   \centering
   \includegraphics[width=7in]{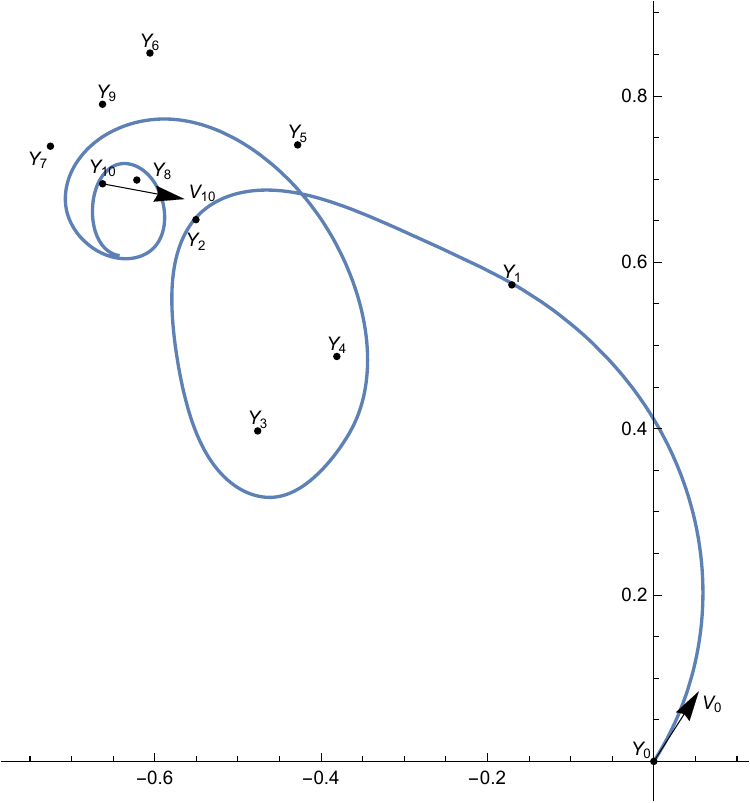} 
   \caption{$y_{\hat \theta}$ for Example \ref{ex3}}
   \label{fig4}
\end{figure}
We see that $y_{\hat \theta}$ is very far from interpolating, except at $T_0,T_1,T_2$. Even though the curve nearly passes through the terminal point $Y_{10}$, it arrives too soon and overshoots. This poor performance can be corrected by more refined sampling or, as illustrated in Example \ref{ex6}, 
by the more sophisticated method of \S \ref{gapsec} which starts with $\hat \theta$. \qed
\end{example}
\section{Some Asymptotic Lemmas}\label{asympsec}
Recall $\theta _j(s)=a_j+b_js+c_js^2+d_js^3$ where $s\in [0,L_j]$ for $1\leq j\leq n$, and $Y,V,T$ are given with $(Y,V)$ is admissible with respect to $T$. Using Mathematica
 to calculate the Taylor polynomial $\hat F_j(s)$ of degree $4$ of $F_j(s):=\int _0^s(\cos u,\sin u)~du)$ about $s=0$, we find 
\begin{lemma}\label{lemm1}  For $1\leq j\leq n$, 
$$\hat F_j(L_j)=L_j\left[ \begin{array}{cc} 
\cos a_j &-\sin a_j\\
\sin a_j& ~\cos a_j 
\end{array}\right]   \left[ \begin{array}{c}\alpha _j \\
\beta _j 
\end{array}\right] $$ 
where 
$$\alpha _j:=1-\frac{b_j^2L_j^2}{6}-\frac{b_jc_jL_j^3}{4}+\frac{(b_j^4-12c_j^2-24b_jd_j)L_j^4}{120}\hbox{~and~}
\beta _j:=\frac{b_jL_j}{2}+\frac{c_jL_j^2}{3}-\frac{(b_j^3-6d_j)L_j^3}{24}-\frac{b_j^2c_jL_j^4}{10}.$$
\qed 
\end{lemma}

\spp
From Lemma \ref{lemm1},  
$$L_jr_j\left[ \begin{array}{c} \cos \omega _j\\
\sin \omega _j
\end{array}\right] ~=~L_jq_j ~=~Y_j-Y_{j-1}~=~F_j(L_j)~=~L_j\left[ \begin{array}{cc} \cos a_j &-\sin a_j\\
\sin a_j& ~\cos a_j 
\end{array}\right]   \left[ \begin{array}{c}\alpha _j \\
\beta _j 
\end{array}\right] +O(\epsilon ^5)~\Longrightarrow $$
\begin{eqnarray}\label{rjeq}r_j^2&=&\alpha _j^2+\beta _j^2+O(\epsilon ^4)\\
\label{omjeq} (\cos (a_j-\omega _j),\sin (a_j-\omega _j))&=&\frac{(\alpha _j,\beta _j)}{\Vert (\alpha _j,\beta _j)\Vert}+O(\epsilon ^5). \end{eqnarray} 
\begin{lemma}\label{lemm1p1}  
\begin{eqnarray}\label{oldeq11}(b_j+c_jL_j+\frac{9}{10}d_jL_j^2)^2+\frac{c_j^2L_j^2}{15}&=&k_j^2(1-\frac{k_j^2L_j^2}{20})+O(\epsilon ^3)\quad \quad \hbox{ for }1\leq j\leq n,\\
\label{ajeq}
a_j+\frac{b_jL_j}{2}+\frac{c_jL_j^2}{3}+\frac{d_jL_j^3}{4}&=&\omega _j+O(\epsilon ^4)\quad \quad \quad \quad \quad \quad \quad \hbox{for }1\leq j\leq n,\\
\label{oldeq10} \frac{b_{j+1}L_{j+1}+b_jL_j}{2}+\frac{c_{j+1}L_{j+1}^2+2c_jL_j^2}{3}+\frac{d_{j+1}L_{j+1}^3+3d_jL_j^3}{4}&=&\omega _{j+1}-\omega _j+O(\epsilon ^4)\quad \quad \quad \quad \hbox{for }1\leq j\leq n-1.
\end{eqnarray}
\end{lemma}

\spp
{\bf Proof:} Using Mathematica to substitute for $r_j^2$ from (\ref{rjeq}) in $k_j^2$ and then expand in powers of $L_j$, we obtain  
\begin{equation}\label{oldeqbefore11}
(b_j+c_jL_j+\frac{9}{10}d_jL_j^2)^2+(\frac{b_j^4}{20}+\frac{c_j^2}{15})L_j^2~=~k_j^2.\end{equation}
where $1\leq j\leq n$. In particular, $b_j^2=k_j^2+O(\epsilon )$. Substituting  $b_j^4=k_j^4+O(\epsilon )$ back in (\ref{oldeqbefore11}) gives (\ref{oldeq11}). 
For (\ref{ajeq}) use Mathematica to Taylor expand the argument of the right hand side of (\ref{omjeq}) in powers of $L_j$. For (\ref{oldeq10}) substitute  in (\ref{eqct})  for $a_j,a_{j+1}$ from (\ref{ajeq}). 
\qed 
\section{Estimating the $b_j$ to $O(\epsilon ^2)$}\label{newnewsec} 
From (\ref{oldeq10}), for $1\leq j\leq n-1$,   
\begin{equation}\label{oldeq10cut}\displaystyle{\frac{b_{j+1}L_{j+1}+b_jL_j}{2}+\frac{c_{j+1}L_{j+1}^2+2c_jL_j^2}{3}=\omega _{j+1}-\omega _j+O(\epsilon ^3)}.\end{equation}
From (\ref{ajeq}) with $j=1$, and from (\ref{aux}), 
\begin{equation}\label{firstbeq}
c_1L_1~=~\frac{3(\omega _1-\nu _0)}{L_1}-\frac{3b_1}{2}+O(\epsilon ^2).
\end{equation}
From (\ref{eqc1}) for $1\leq j\leq n-1$, 
$$c_jL_j=\frac{b_{j+1}-b_j}{2}+O(\epsilon ^2).$$ 
From (\ref{ajeq}) with $j=n$, and from (\ref{aux}),  
\begin{equation}\label{cnest}
\omega _n-\frac{b_nL_n}{2}-\frac{c_nL_n^2}{3}~=~a_n+O(\epsilon ^3)~=~\nu _n-b_nL_n-c_nL_n^2~\Longrightarrow~c_nL_n~=~-\frac{3(\omega _n-\nu _n)}{2L_n}-\frac{3b_n}{4} +O(\epsilon ^2).
\end{equation}
Substituting for $c_1L_1$ in (\ref{firstbeq}), and for $c_1L_1,c_2L_2,\ldots ,c_nL_n$ in (\ref{oldeq10cut}) with $1\leq j\leq n-1$,   we obtain the system 
\begin{eqnarray*}
2L_1b_1+L_1b_{2}&=&6(\omega _1-\nu _0)+O(\epsilon ^3),\\
L_{j-1}b_{j-1}+2(L_{j-1}+L_j)b_j+L_{j+1}b_{j+1}&=&6(\omega _j-\omega _{j-1})+O(\epsilon ^3)\quad \hbox{for }2\leq j\leq n-2,\\
L_{n-1}b_{n-1}+(2L_{n-1}+\frac{3L_n}{2})b_n&=&3(3\omega _n-\nu _n-2\omega _{n-1})+O(\epsilon ^3),
\end{eqnarray*}
of $n$ linear equations for $b:=(b_1,b_2,\ldots ,b_n)$, namely 
\begin{equation}\label{btilapp} \tilde {\bf T}~b~=~\tilde {\bf R}+O(\epsilon ^3)~=~\tilde {\bf T}~\tilde b+O(\epsilon ^3)\end{equation}
where $\tilde {\bf T}$ is the $n\times n$ tridiagonal matrix and $\tilde {\bf R}, \tilde b$ are the $n$-dimensional vectors all defined in \S \ref{newsec}. 

\begin{lemma}\label{btillem}
$b=\tilde b+O(\epsilon ^2)$.
\end{lemma}

\spp
{\bf Proof:} Let $D_{\tilde {\bf T}}$ be the diagonal $n\times n$ matrix with the same diagonal entries as $\tilde {\bf T}$. Because $A_m\epsilon <L_j<A_M\epsilon $, we see that $\Vert D_{\tilde {\bf T}}^{-1}\Vert _{\infty}\leq 1/(2A_m\epsilon )$. By (\ref{btilapp}),  $\tilde {\bf T}~(b-\tilde b)=O(\epsilon ^3)$, and therefore  
 $$\Vert D_{\tilde {\bf T}}^{-1}\tilde {\bf T} ~(b-\tilde b)\Vert ~=~O(\epsilon ^2).$$
Because $\tilde {\bf T}$ is strictly diagonally dominant by rows with dominance factor $1/2$, 
$D_{\tilde {\bf T}}^{-1}\tilde {\bf T}$ has condition number $\leq 3$. 
\qed 

\section{Proof of (\ref{bjhat})}\label{newnewsec1}
By (\ref{eqc1}) and Lemma \ref{btillem}, 
$$c_jL_j~=~\frac{b_{j+1}-b_j}{2}+O(\epsilon ^2)~=~\frac{\tilde b_{j+1}-\tilde b_j}{2}+O(\epsilon ^2)$$ 
for $1\leq j\leq n-1$.  Similarly by (\ref{cnest}) and Lemma \ref{btillem}, ~$\displaystyle{c_nL_n=\frac{3(\omega _n-\nu _n)}{2L_n}-\frac{3\tilde b_n}{4}+O(\epsilon ^2)}$.~ Substituting for the $c_jL_j$ in (\ref{oldeq11}),  
\begin{eqnarray*}
(b_j+c_jL_j+\frac{9}{10}d_jL_j^2)^2&=&k_j^2(1-\frac{k_j^2L_j^2}{20})-\frac{(\tilde b_{j+1}-\tilde b_j)^2}{60}+O(\epsilon ^3)\quad \hbox{ for }1\leq j\leq n-1, \hbox{ and}\\
(b_n+c_nL_n+\frac{9}{10}d_nL_n^2)^2&=&k_n^2(1-\frac{k_n^2L_n^2}{20})-\frac{9}{60}\left( \frac{\omega _n-\nu _n}{L_n}+\frac{\tilde b_n}{2}\right) ^2+O(\epsilon ^3),
\end{eqnarray*}
where the right hand sides are known to $O(\epsilon ^3)$. 
By Lemma \ref{bdsklem}, for $\epsilon $ small the $k_j$ are bounded away from $0$ for all $1\leq j\leq n$. 
Therefore
$$b_j+c_jL_j+\frac{9}{10}d_jL_j^2~=~\sigma _j\vert \rho _j\vert +O(\epsilon ^3)~=~\sigma _jk_j+O(\epsilon ^2)$$
where $\sigma _j=\pm 1$, and the $\rho _j$ are defined as in (\ref{rhojdf}), (\ref{rhondf}). In particular, by Lemma \ref{bdsklem} $\sigma _j={\rm sign} (b_j)$ for $1\leq j\leq n$. By (\ref{eqc1}), $b_{j}=b_{j-1}+O(\epsilon )$ 
for $2\leq j\leq n$. Because the $k_j$ are bounded away from $0$ so are the $\vert b_j\vert$.  So the $\sigma _j$ are all equal.  Recall the definition of $\sigma$ in \S \ref{newsec}, namely $\sigma := {\rm sign} (\det [ V_0~\vdots ~Y_1-Y_0])$. 

\spp
\begin{lemma}\label{siglem} For $\epsilon >0$ small, we have $\sigma _j=\sigma$ where $1\leq j\leq n$.  
\end{lemma}

\spp
{\bf Proof:} By (\ref{eqc1}), $b_1=b_2+O(\epsilon)$. Then, from the first equation in system (\ref{Ttileq}), 
$$3b_1~=~2b_1+b_2+O(\epsilon )~=~2\tilde b_1+\tilde b_2+O(\epsilon )~=~\frac{6(\omega _1-\nu _0)}{L_1}+O(\epsilon )~\Longrightarrow ~\sigma _1~=~{\rm sign}(\omega _1-\nu _0),$$
where we use $\sigma _1={\rm sign }(b_1)$.  
Now $V_0=x^{(1)}(T_0)$ and $Y_j=x(T_j)$ for $0\leq j\leq n$, where $x:[T_0,T_n]\rightarrow \R ^2$ is a convex generator. Write $x^{(1)}(t)=(\cos \psi (t),\sin \psi (t))$ 
where $\psi $ is $C^1$ and $\psi (T_0)=\nu _0$. By convexity, 
$\vert \psi ^{(1)}(T_0)\vert =\Vert x^{(2)}(T_0)\Vert >C>0$. Taylor expanding $x$ to second order in $t$ about $T_0$, and recalling the definition of $\omega _j$, 
$$L_1r_1(\cos \omega _1,\sin \omega _1)~=~Y_1-Y_0~=~L_1(\cos \nu _0,\sin \nu _0)+\frac{L_1^2}{2}\psi ^{(1)}(T_0)(-\sin  \nu _0,\cos \nu _0)+O(\epsilon ^3)~\Longrightarrow$$
$$L_1r_1\sin (\omega _1-\nu _0)~=~\frac{L_1^2}{2}\psi ^{(1)}(T_0)+O(\epsilon ^3)~=~\det [V_0~\vdots ~Y_1-Y_0]+O(\epsilon ^3)~\Longrightarrow~ {\rm sign}(\omega _1-\nu _0)~=~\sigma .$$
So $\sigma _n=\sigma _{n-1}=\ldots =\sigma _1=\sigma$. $\qed$

\spp
From Lemma \ref{siglem} we now have, for $1\leq j\leq n$, 
\begin{equation}\label{rheq}
b_j+c_jL_j+\frac{9}{10}d_jL_j^2~=~\rho _j +O(\epsilon ^3).
\end{equation}

\spp
By (\ref{eqc1}), (\ref{rheq}), for each $1\leq j\leq n-1$, 
\begin{eqnarray*}
2c_jL_j+3d_jL_j^2&=&b_{j+1}-b_j+O(\epsilon ^3),\\
c_jL_j+\frac{9}{10}d_jL_j^2&=&\rho _j -b_j+O(\epsilon ^3),
\end{eqnarray*} 
where the right hand sides are presently known only with $O(\epsilon ^2)$ errors.  Nevertheless, solving these systems with $1\leq j\leq n-1$, 
\begin{eqnarray}
\label{acjeq}c_jL_j&=&\frac{-7b_j-3b_{j+1}}{4}~+\frac{5\rho _j}{2}+O(\epsilon ^3),\\
\label{adjeq}d_jL_j^2&=&\frac{5(b_j+b_{j+1})}{6}\quad -\frac{5\rho _j}{3}+O(\epsilon ^3).
\end{eqnarray}
By (\ref{ajeq}) with $j=n$ and (\ref{aux}), and by (\ref{rheq}) with $j=n$, 
\begin{eqnarray*}8c_nL_n+9d_nL_n^2&=&\frac{12(\nu _n-\omega _n)}{L_n}-6b_n+O(\epsilon ^3),\\
10c_nL_n+9d_nL_n^2&=&10\rho _n-10b_n+O(\epsilon ^3).
\end{eqnarray*}
Solving, 
\begin{eqnarray}
\label{acneq}c_nL_n&=& \frac{6(\omega _n-\nu _n)}{L_n}-2b_n+5\rho _n+O(\epsilon ^3),\\
\label{adneq}d_nL_n^2&=&-\frac{20(\omega _n-\nu _n)}{3L_n}+\frac{10b_n}{9}-\frac{40\rho _n}{9}+O(\epsilon ^3)  .
\end{eqnarray}
In particular, substituting for $c_1L_1,d_1L_1^2$ in (\ref{ajeq}) with $j=1$, 
$$3L_1b_1-L_1b_2~=~24(\omega _1-\nu _0)-10 L_1\rho _1+O(\epsilon ^4).$$
\spp
Substituting for $c_jL_j,d_jL_j^2$ where $1\leq j\leq n$ in (\ref{oldeq10}), gives $n-1$ linear equations for the remaining components of $b$, namely
\begin{eqnarray*}
-L_jb_j+3(L_j+L_{j+1})b_{j+1}-L_{j+1}b_{j+2}&=&24(\omega _{j+1}-\omega _j)-10(L_j\rho _j+L_{j+1}\rho _{j+1})+O(\epsilon ^4)\quad \hbox{for }1\leq j\leq n-2, \quad \hbox{and}\\
-3L_{n-1}b_{n-1}+(9L_{n-1}+8L_n)b_n&=&24(2\omega _n+\nu _n-3\omega _{n-1})-10(3L_{n-1}\rho _{n-1}+4L_n\rho _n)+O(\epsilon ^4).
\end{eqnarray*} 

\spp
So $\hat {\bf T}~b=\hat R +O(\epsilon ^4)$ and  (\ref{bjhat}) follows, in the same way as for Lemma \ref{btillem}.  $\qed$ 
\section{Proof of Theorem \ref{thm0}}\label{completepfsec}
As noted at the end of \S \ref{newsec}, now that (\ref{bjhat}) is proved it suffices to verify (\ref{ajhat}), (\ref{cjhat}), (\ref{djhat}). 
Using (\ref{bjhat}) in (\ref{acjeq}), (\ref{adjeq}), then in (\ref{acneq}), (\ref{adneq}),  
\begin{eqnarray*}
c_jL_j&=&\frac{-7b_j-3b_{j+1}}{4}~+\frac{5\rho _j}{2}+O(\epsilon ^3)~=~\frac{-7\hat b_j-3\hat b_{j+1}}{4}~+\frac{5\rho _j}{2}+O(\epsilon ^3)~=~\hat c_jL_j+O(\epsilon ^3)\quad \hbox{for }1\leq j\leq n-1,\\
\label{adjeq}d_jL_j^2&=&\frac{5(b_j+b_{j+1})}{6}\quad -\frac{5\rho _j}{3}+O(\epsilon ^3)~=~
\frac{5(\hat b_j+\hat b_{j+1})}{6}\quad -\frac{5\rho _j}{3}+O(\epsilon ^3)~=~\hat d_jL_j^2+O(\epsilon ^3)\quad \hbox{for }1\leq j\leq n-1,\\
c_nL_n&=& \frac{6(\omega _n-\nu _n)}{L_n}-2b_n+5\rho _n+O(\epsilon ^3)~=~\frac{6(\omega _n-\nu _n)}{L_n}-2\hat b_n+5\rho _n+O(\epsilon ^3)~=~ \hat c_nL_n+O(\epsilon ^3),\\
\label{adneq}d_nL_n^2&=&-\frac{20(\omega _n-\nu _n)}{3L_n}+\frac{10b_n}{9}-\frac{40\rho _n}{9}+O(\epsilon ^3)~=~-\frac{20(\omega _n-\nu _n)}{3L_n}+\frac{10\hat b_n}{9}-\frac{40\rho _n}{9}+O(\epsilon ^3)~=~\hat d_nL_n^2+O(\epsilon ^3),  
\end{eqnarray*}
according to definitions (\ref{bcjeq}), (\ref{bdjeq}),  (\ref{bcneq}), (\ref{bdneq}) in \S \ref{newsec}. This proves (\ref{cjhat}), (\ref{djhat}). Finally by (\ref{aux}) $a_1=\nu _0=\hat a_1$, and by (\ref{eqct}) $a_{j+1}=a_j+b_jL_j+c_jL_j^2+d_jL_j^3=
\hat a_j+\hat b_jL_j+\hat c_jL_j^2+\hat d_jL_j^3+O(\epsilon ^4)=\hat a_{j+1}+O(\epsilon ^4)$ according to (\ref{bajeq}) for $1\leq j\leq n-1$. So (\ref{ajhat}) is proved by induction. $\qed$

\section{Closing Gaps}\label{gapsec}
As noted in \S \ref{appsec}, the $C^2$ unit-speed curve $y_{\hat \theta}$ does not interpolate $(Y,V)$ {\em exactly} at 
$T$. Usually there are $O(\epsilon ^5)$ errors in the interpolation conditions. The left-interpolant $\bar y_{\hat \theta}$ is not much better, except to illustrate the magnitudes and locations of errors, because $\bar y_{\hat \theta}$ is not even $C^0$. 

\spp
When sampling is sufficiently fine, namely for sufficiently small $\epsilon >0$, the curves 
$y_{\hat \theta}$ and $\bar y_{\hat \theta}$ are almost indistinguishable, and either of these  unit-speed approximate interpolants might be used in practice. This is already true for some 
parts of $y_{\hat \theta}$ in Example \ref{ex1} and for most of $y_{\hat \theta}$ in Example \ref{ex2}. However $y_{\hat \theta}$ is seriously problematic in Example \ref{ex3}. 

\spp
If finer sampling is inconvenient or impossible then, after finding $\hat \theta$, there is a second calculation we can do to correct the kinds of errors observed in Example \ref{ex3}, and also to a smaller extent in Examples \ref{ex1}, \ref{ex2}. 

\spp
We rewrite the conditions for a second order spiral spline interpolant as a nonlinear system of equations, then solve numerically using $\hat \theta $ to generate an initial guess. This is easier than finding $\hat \theta$ and simpler to code, but not so remarkably fast. It is unlikely to fail except when sampling is so coarse that $\hat \theta$ does not give a good initial guess. Even in the case of Example \ref{ex3} the rough estimate $\hat \theta$ is enough to get started. 

\spp
Given $(Y,V)$ admissible with respect to $T$, first calculate $\hat \theta$ as in \S \ref{newsec}.  
Then, for $u,v\in \R ^n$, define $a_j,b_j,c_j,d_j\in \R$ for $1\leq j\leq n$ by $c_j:=u_j$ and $d_j:=v_j$ 
for all $1\leq j\leq n$, and $a_1:=\nu _0$. Using (\ref{eqct}) and (\ref{eqc1}) $j-1$ times where $1\leq j\leq n$, 
\begin{eqnarray}\label{bjb1eq}
b_j&=&b_1+2\sum_{k=1}^{j-1}u_kL_k+3\sum_{k=1}^{j-1}v_kL_k^2\quad \hbox{and}\\
\label{ajbkeq}a_j&=&\nu _0+\sum_{k=1}^{j-1}b_kL_k +\sum_{k=1}^{j-1}u_kL_k^2+\sum_{k=1}^{j-1}v_kL_k^3.\end{eqnarray}
Using (\ref{bjb1eq}) to substitute for $b_k$ in (\ref{ajbkeq}),
\begin{equation}
\label{ajeqf}a_j~=~\nu _0+(T_{j-1}-T_0)b_1+2\sum_{k=1}^{j-1}\sum_{i=1}^{k-1}u_iL_iL_k+3\sum_{k=1}^{j-1}\sum_{i=1}^{k-1}v_iL_i^2L_k +\sum_{k=1}^{j-1}u_kL_k^2+\sum_{k=1}^{j-1}v_kL_k^3
\end{equation}
where $2\leq j\leq n$. Similarly, using the second part of (\ref{aux}) in place of (\ref{eqct}), 
\begin{eqnarray}
\nonumber \nu _n&=&\nu _0+(T_{n}-T_0)b_1+2\sum_{k=1}^{n}\sum_{i=1}^{k-1}u_iL_iL_k+3\sum_{k=1}^{n}\sum_{i=1}^{k-1}v_iL_i^2L_k +\sum_{k=1}^{n}u_kL_k^2+\sum_{k=1}^{n}v_kL_k^3 ~\Longrightarrow \\
\label{b1eq}b_1&=&\frac{ \nu _n-\nu _0-2\sum_{k=1}^{n}\sum_{i=1}^{k-1}u_iL_iL_k-3\sum_{k=1}^{n}\sum_{i=1}^{k-1}v_iL_i^2L_k -\sum_{k=1}^{n}u_kL_k^2-\sum_{k=1}^{n}v_kL_k^3}{T_n-T_0}.
\end{eqnarray}
Then $a_j$ and $b_j$ are also determined for $2\leq j\leq n$, by substitution for $b_1$ in (\ref{ajeqf}) 
and in (\ref{bjb1eq}). For $1\leq j\leq n$ define $\theta _j:[0,L_j]\rightarrow \R$ by $\theta _j(t):=a_j+b_jt+c_jt^2+d_jt^3$. Define $\theta :[T_0,T_n]\rightarrow \R$ by $\theta (T_0):=\nu _0$ 
and $\theta (t):=\theta _j(t-T_{j-1})$ for $t\in (T_{j-1},T_j]$. 
Because the $a_j,b_j,c_j,d_j$ satisfy (\ref{eqct}), (\ref{eqc1}), as well as the auxiliary conditions (\ref{aux}), the following lemma is easily verified. 
\begin{lemma}\label{uvlem}   $\theta :[T_0,T_n]\rightarrow \R$ is a cubic polynomial spline with  
$\theta (T_0)=\nu _0$ and $\theta (T_n)=\nu _n$. \qed
\end{lemma}

\spp
So for any $u,v\in \R ^n$,  $y_\theta :[T_0,T_n]\rightarrow \R ^2$ is a $C^2$ second order spiral spline 
satisfying $y_\theta ^{(1)}(T_0)=V_0$ and $y_\theta ^{(1)}(T_n)=V_n$.
For $1\leq j\leq n$ define  
$z_j:\R ^n\times \R ^n\rightarrow \R ^2$ by 
\begin{equation}\label{zsteq}\displaystyle{z_j(u,v):=\int _0^{L_j}(\cos \theta _j(t),\sin \theta _j(t))~dt},\end{equation}
and define $z:\R ^n\times \R ^n\rightarrow (\R ^2)^n$ by $z(u,v):=(z_1(u,v),z_2(u,v),\ldots ,z_n(u,v))$. 

\spp
Then $y_\theta$ interpolates $(Y,V)$ at $T$ precisely when $z_j(u,v)=Y_j-Y_{j-1}$ 
for all $1\leq j\leq n$, namely when $(u,v)\in \R ^n\times \R ^n$ satisfies the nonlinear system of $2n$ scalar equations in $2n$ scalar unknowns
\begin{equation}\label{zeq} z(u,v)=z^*:=(Y_1-Y_0,Y_2-Y_1,\ldots ,Y_n-Y_{n-1})\in (\R ^2)^n\end{equation}
whose right hand sides are found from $Y$. The key to solving such a system is a satisfactory initial guess $(\hat u,\hat v)$ to a solution $(u,v)$. Set $\hat u=\hat c$ and $\hat v=\hat d$, 
where the $\hat a_j,\hat b_j,\hat c_j,\hat d_j$ are found in \S \ref{newsec}. By Corollary \ref{cor0}, $z(\hat u,\hat v)=z^*+O(\epsilon ^5)$.

\spp
In practice the integrals on the right hand side of (\ref{zsteq}) are not easy to express as functions of $u,v$. So for numerical computations we replace the $z_j$ by approximations using the composite Simpson's Rule. Then $z$ becomes an explicit function of $(u,v)$ and (\ref{zeq}) is solved using Mathematica's FindRoot with $(\hat u,\hat v)$ as an initial guess.
%

\spp
\begin{example}\label{ex4} Using composite Simpson with $6$ intervals to estimate the $z_j$, FindRoot 
took around $0.3$ seconds to obtain a numerical solution $\theta$ of (\ref{zeq}) for the data in Example \ref{ex1}. Figure \ref{fig5} shows the interpolant $y_\theta$ (yellow), together with the initial estimate $y_{\hat \theta}$ and data from Figure \ref{fig2}. 
\begin{figure}[htbp] 
   \centering
   \includegraphics[width=7in]{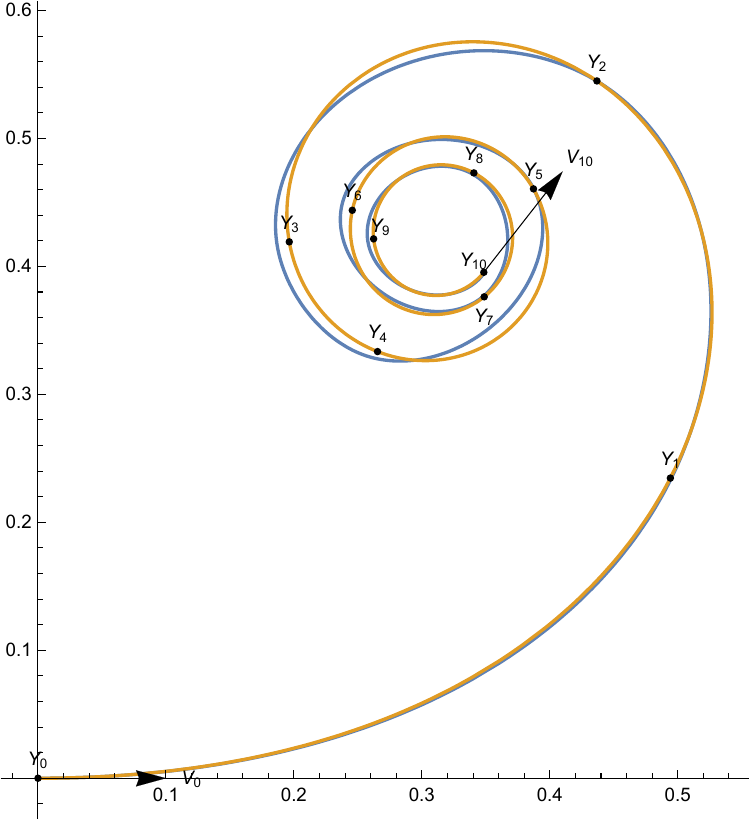} 
   \caption{$y_{ \theta}$ (yellow) and  $y_{\hat \theta}$  for Example \ref{ex4}}
   \label{fig5}
\end{figure}
\qed
\end{example}

\begin{example}\label{ex5} Perhaps it is no surprise that $y_\theta$ is successful for the data of  Example \ref{ex2}, because  $y_{\hat \theta}$ was already nearly interpolating.  Figure \ref{fig6} shows $y_\theta$ (yellow), together with $y_{\hat \theta}$ and the data. FindRoot took some 0.3 seconds to find $y_\theta$. 
\begin{figure}[htbp] 
   \centering
   \includegraphics[width=7in]{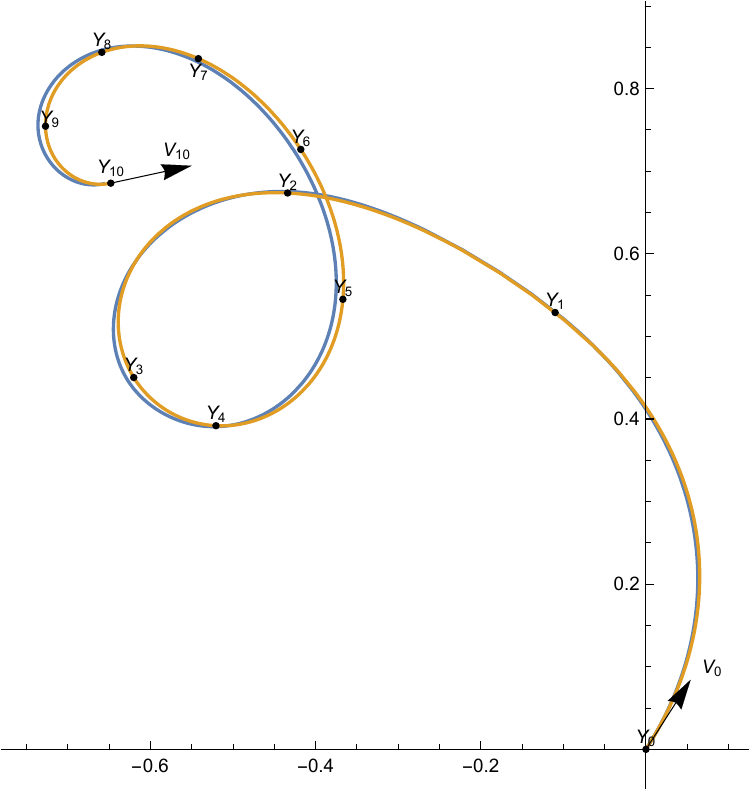} 
   \caption{$y_{ \theta}$ (yellow) and $y_{\hat \theta}$ in Example \ref{ex5}}
   \label{fig6}
\end{figure}
\qed
\end{example}

\begin{example}\label{ex6}{ex6}
More impressively, it takes some 0.3 seconds to improve $y_{\hat \theta}$ in Example \ref{ex3} to the interpolant $y_{\theta}$ shown (yellow) in Figure \ref{fig7}. Whereas $y_{\hat \theta}$ comprehensively failed to interpolate, the $C^1$ cubic spline $\hat \theta :[T_0,T_n]\rightarrow \R$ retains considerable informative power. This power, masked to some degree\footnote{Interpolation errors accumulate along the trajectory. The related curve $\bar y_{\hat \theta}$ of \S \ref{appsec} avoids this defect, and is better for diagnostics, but fails to be continuous at interpolation points.} by the second order spiral spline $y_{\hat \theta}:[T_0,T_n]\rightarrow \R ^2$, is key to the method of the present section \S \ref{gapsec}.  
\begin{figure}[htbp] 
   \centering
   \includegraphics[width=7in]{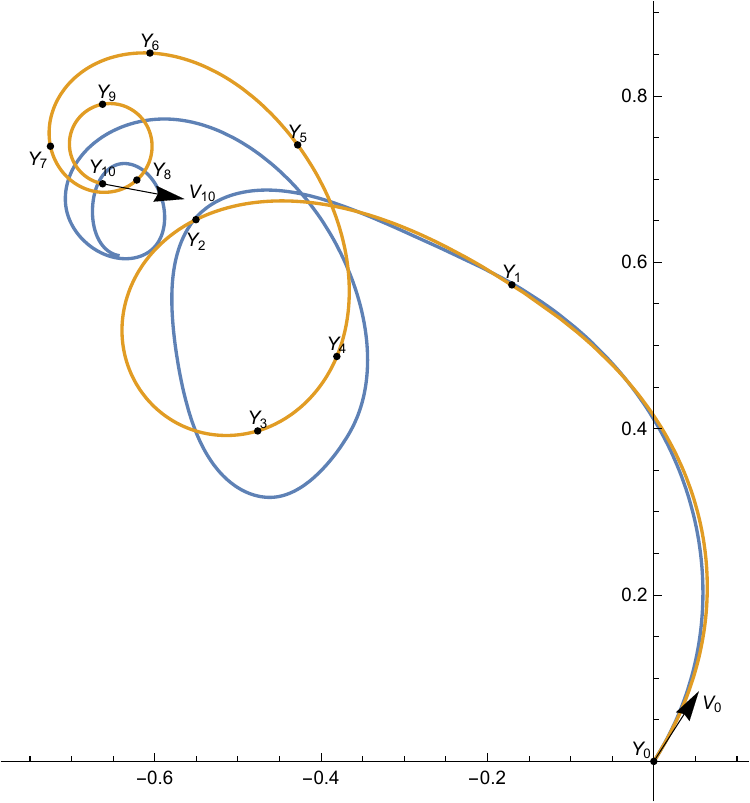} 
   \caption{$y_{ \theta}$ (yellow) and $y_{\hat \theta}$ for Example \ref{ex6}}
   \label{fig7}
\end{figure}
\qed
\end{example}

\section*{Conclusion}
Given interpolation data $(Y,V)$ for convex unit-speed $C^2$ curves in $\R ^2$, relative to some sequence $T$, we develop a method for interpolating $(Y,V)$ at $T$ by a $C^2$ second order spiral spline $y_\theta :[T_0,T_n]\rightarrow \R ^2$. The method is developed using asymptotic arguments to define a pair of tridiagonal linear systems whose solutions define a real-valued $C^1$ cubic polynomial spline $\hat \theta $ with prescribed values at $T_0,T_n$. For sufficiently finely sampled data, $\hat \theta$ already defines a $C^2$ second order spiral spline $y_{\hat \theta}$ that approximately interpolates $(Y,V)$ at $T$. More importantly, $\hat \theta$ is used to start a numerical method that finds another $C^1$ cubic polynomial spline $\theta$. Then $y_\theta$ is the desired $C^2$ second order interpolating spiral spline.



\bibliographystyle{elsarticle-num}
\bibliography{<your-bib-database>}

\begin{thebibliography}{00}


\bibitem{dan}Di Antonio L. The fabric of the universe is most perfect: Euler's research on elastic curves, in Euler at, 300: an appreciation, 239--260, Mathematical Association of America (2007). 


\bibitem{audoly} Audoly, B., Pomeau, Yves, Elasticity and Geometry: From Hair Curls to the Non-linear Response of Shells, Oxford University Press (2010). 

\bibitem{baran} Baran, Ilya, Lehtinen, Jaako and Popovic, Jovan, Sketching Clothoid Splines Using Shortest Paths, Eurographics 0 (1981).

\bibitem{bertails} Florence Bertails-Descoubes, Super-Clothoids, Eurographics 31 (2012) 

\bibitem{birkhoff} G. Birkhoff, H. Burchard and D. Thomas, Nonlinear Interpolation by Splines, Pseudosplines and Elastica, General Motors Research Laboratories Report 468, (1965), 1--13.

\bibitem{bruckstein} Bruckstein, Alfred M., Holta Robert, J., and Netravalia, Arun, N., 
Discrete elastica,  Applicable Analysis 78 (2001) 453--485.

\bibitem{brunnett} Brunnett, Guido and Wendt, J\"org, A Univariate Method for Plane Elastic Curves, 
Computer Aided Geometric Design 14 (1997) 273--292. 

\bibitem{coope} Coope, Ian, D., Curve Fitting With Nonlinear Spiral Splines, Department of Mathematics, University of Canterbury NZ, No 63 (1991), 1--14.

\bibitem{dai} Dai,  Ran and Cochran,  John E. Jr. (2010) Path Planning and State Estimation for Unmanned
Aerial Vehicles in Hostile Environments,
Journal of Guidance,  Control  and Dynamics 33 (2), 595--601.


\bibitem{deBoor} de Boor, Carl, A Practical Guide to Splines, Applied mathematical Sciences 27, Springer (1978). 

\bibitem{eberly} Eberly, David, Moving Along a Curve with Specified Speed, {\tt http://www.geometrictools.com/} (2007) 1--15.

\bibitem{edwards} Edwards, J.A., Exact equations of the nonlinear spline, 
ACM Transactions on Mathematical Software 18 (1992), 174--192.


\bibitem{jerome} Golomb, Michael and Jerome, Joseph, Equilibria of the Curvature Functional and Manifolds of Interpolating Spline Curves, SIAM J. on Mathematical Analysis 13 (3) (1982), 412--458.


\bibitem{Kelly} Alonzo Kelly and Brian Nagy, Reactive Nonholonomic Trajectory Generation via Parametric Optimal Control, International Journal of Robotics Research 22 (2003), 583--601. 
  
\bibitem{levien} Levien, Raph, The Elastica: a Mathematical History, Technical Report 
No. UCB/EECS-2008-103\\ 
{\tt http://www.eecs.berkeley.edu/Pubs/TechRpts/2008/EECS-2008-103.html} 
August 23, 2008. 

\bibitem{leefor} E.H. Lee and G.E. Forsythe, Variational Study of Nonlinear Spline Curves, SIAM Review 15 (1973), 120--133.

\bibitem{linners} Linners, Anders, Unified Representations of Nonlinear Splines, Journal of Approximation Theory 84 (1996), 315--350.  

\bibitem{looker} Looker, Jason R., Constant Speed Interpolating Paths, DSTO Defence Science and Technology Organisation, AR 014-939. DSTO-TN-0989, March 2011. 

\bibitem{meek0} D.S. Meek and R.S.D. Thomas, A Guided Clothoid Spline, Computer Aided Geometric Design 8 (1991) 163--174. 

\bibitem{meek} D.S. Meek and D.J. Walton, An Arc Spline Approximation to a Clothoid, Journal of Computational and Applied Mathematics 170 (2004), 59--77. 



\bibitem{peterson} Peterson, John W., Arc-Length Parameterization of Spline Curves {\tt http://www.saccade.com/writing/graphics/RE-PARAM.PDF}


\bibitem{singer} Singer, David A., Lectures on Elastic Curves and Rods, in AIP Conference Proceedings 1002, Curvature and Variational Modelling in Physics and Biophysics, Santiago de Compostela, Spain, 17--18 September 2007 (eds \'Oscar J. Garay, Eduardo Garcia-R\'io, Ram\'on V\'azquez-Lorenzo), 3--32.

\bibitem{stoer} Stoer, Josef, Curve Fitting With Clothoidal Splines, J. Research of the National Bureau of Standards 87 No. 4 (1982) 318--346. 















\bibitem{WKA} Arc-Length Parameterized Spline Curves for Real-Time Simulation, Wang, Hongling, Keaney, Joseph, Atkinson, Kendall, Curve and Surface Design: Saint Malo (2002) Lyche, T., Mazure, M.-L., Schumacher, L. (eds) 387--396. 



\end{thebibliography}



\end{document}